\documentclass[a4paper,12pt,leqno]{amsart}

\usepackage{amssymb,latexsym}
\usepackage{mltex}
\usepackage{amsmath}
\usepackage{epsfig}
\usepackage{graphics}
\usepackage{pb-diagram}
 \usepackage{lscape} 
\usepackage{enumerate}
\usepackage{amsthm}

\newtheorem{thm}{Theorem}     
\newtheorem{cor}{Corollary}     
\newtheorem{lem}{Lemma}[section]
\newtheorem{prop}[lem]{Proposition}

\newcommand{\Mndrk}{\mathcal{M}(n,\delta,R,k)}
\newcommand{\nablab}{\overline{\nabla}}
\newcommand{\expo}{\exp_{p_0}}
\newcommand{\Ric}{\mathrm{Ric}}
\newcommand{\diff}{\mathrm{d}\,}
\newcommand{\td}{t_{\delta}}
\newcommand{\sd}{s_{\delta}}
\newcommand{\cd}{c_{\delta}}
\newcommand{\lgra}{\longrightarrow}
\newcommand{\Hinf}{||H||_{\infty}}
\newcommand{\Hkinf}{||H_k||_{\infty}}
\newcommand{\Hkp}{||H_k||_{2p}}
\newcommand{\HHkp}{||H_k||^2_{2p}}

\newcommand{\Binf}{||B||_{\infty}}

\newcommand{\iid}{\mathrm{Id}\,}

\newcommand{\ddiv}{\mathrm{div}\,}

\newcommand{\iinf}{\mathrm{inf}\,}

\newcommand{\sscal}{\mathrm{Scal}\,}

\newcommand{\trace}{\mathrm{tr\,}}

\newcommand{\lto}{\ensuremath{\longrightarrow}}

\newcommand{\Ss}{\mathbb{S}}
\newcommand{\HH}{\mathbb{H}}

\newcommand{\R}{\mathbb{R}}

\newcommand{\Md}{\mathbb{M}^{n+1}(\delta)}

\newcommand{\base}{\{e_1,\ldots ,e_n\}}  
 
\newcommand{\function}[5]
{\begin{eqnarray*}\begin{array}{r@{}ccl}
 #1\;\colon\;  & #2 &\lto & #3 \\[.05cm]  
  & #4 &\longmapsto  & #5 
\end{array}\end{eqnarray*}
}
\newcommand{\beqt}{\begin{equation}}  \newcommand{\eeqt}{\end{equation}}
\newcommand{\bal}{\begin{align}}      \newcommand{\eal}{\end{align}}
\newcommand{\ba}{\begin{array}}      \newcommand{\ea}{\end{array}}
\newcommand{\bc}{\begin{center}}     \newcommand{\ec}{\end{center}}
\newcommand{\be}{\begin{enumerate}}  \newcommand{\ee}{\end{enumerate}}
\newcommand{\beq}{\begin{eqnarray}}  \newcommand{\eeq}{\end{eqnarray}}
\newcommand{\beQ}{\begin{eqnarray*}} \newcommand{\eeQ}{\end{eqnarray*}}
\newcommand{\bi}{\begin{itemize}}    \newcommand{\ei}{\end{itemize}}
\newcommand{\bt}{\begin{tabular}}    \newcommand{\et}{\end{tabular}}
\newcommand{\finpreuve}{\hfill\square\\}
\title[Extrinsic radius pinching]{Extrinsic radius pinching in space forms of nonnegative sectional curavture}
\author{Julien ROTH}

\begin{document}
\maketitle
\begin{center}
Institut \'Elie Cartan, UMR 7502\\
Nancy-Universit\'e, CNRS, INRIA\\
B.P. 239, 54506 Vand\oe uvre l\`es Nancy Cedex, France
\end{center}
\begin{center}
roth@iecn.u-nancy.fr
\end{center}
\begin{abstract}
We first give new estimates for the extrinsic radius of compact hypersurfaces of the Euclidean space $\R^{n+1}$ and the open hemisphere  in terms of high order mean curvatures. Then we prove pinching results corresponding to theses estimates. We show that under a suitable pinching condition, $M$ is diffeomorphic and almost isometric to an $n$-dimensional sphere.\\\\
\end{abstract}

%%%%%%%%%%%%%%%%%%%%%%%%%%%%%%%%%%%%%%%%%%%%%%%%%%%%%%%%%%%%%%%%%%%%%%%%
%%%%%%%%%%%%%%%%%%%%%%%%%%%%%%%%%%%%%%%%%%%%%%%%%%%%%%%%%%%%%%%%%%%%%%%%
\nocite{*}
\section{Introduction}
\label{introduction}
Let $(M^n,g)$ be a compact, connected and oriented $n$-dimensional Riemannian manifold without boundary isometrically immersed by $\phi$ into the $(n+1)$-dimensional simply connected space-form $(\Md,g_{can})$ of sectional curvature $\delta$ with $n\geqslant2$. 
The extrinsic radius of $(M,g)$ is the real number defined as follows
$$R=R(M)=\iinf\left\lbrace r>0\big|\ \exists x\in\Md\ \text{s.t.}\ \phi(M)\subset B(x,r)\right\rbrace,$$
where  $B(x,r)\ \big(${\it resp.} $\overline{B}(x,r)$ and $S(x,r)\big)$ is the open ball $\big(${\it resp.} the closed ball and the sphere$\big)$ of center $x$ and radius $r$ in $\Md$. An immediate consequence of the above definition is that there exists $p_0\in\Md$ such that
$$\phi(M)\subset\overline{B}(p_0,R)\quad\text{and}\quad\phi(M)\cap S(p_0,R)\neq\emptyset.$$
Moreover, the extrinsic radius is bounded from below in terms of the mean curvature. More precisely, we have the following estimate obtained by comparing the shape operators for hypersurfaces that have a contact point and where one is outside the other (see \cite{Ami,BK,HK} for details)
\beqt\label{lowerboundh}
\td(R)\geqslant\frac{1}{\Hinf},
\eeqt
where $\td(s)=\left\lbrace \begin{array}{ll}
 \frac{1}{\sqrt{\delta}}\tan(\sqrt{\delta}s)& \text{if}\ \delta>0 \\ 
 s& \text{if}\ \delta=0 \\ 
 \frac{1}{\sqrt{-\delta}}\tanh(\sqrt{-\delta}s)& \text{if}\ \delta<0
\end{array}\right. $\\
and $H$ is the mean curvature of the immersion.
\indent
Note that if $\delta>0$, the image $\phi(M)$ is assumed to be contained in a ball of radius less than $\frac{\pi}{2\sqrt{\delta}}$ which is equivalent to the fact that $M$ lies in the open hemisphere $\Ss^{n+1}_+(\delta)$.\\
\indent
Since the equality case is characterized by the geodesic spheres of radius $R$, the question of the associated pinching problem was asked, {\it i.e.}, what happens in the case of almost equality. \\
\indent
Many pinching results are known for intrinsic geometric invariants defined on Riemannian manifold with positive Ricci curvature, as the intrinsic diameter, the volume, the radius or the first eigenvalue of the Laplacian (\cite{Ili,Pet,Col1,Col2,Cro,Wu,Aub}).\\
\indent Nethertheless, few results are known about pinching problems in the extrinsic case. In \cite{CG}, B. Colbois and J.F. Grosjean give a first result about the first eigenvalue $\lambda_1(M)$ of the Laplacian. More precisely, they proved that there exists a constant $C$ depending on $n$ and the $L^{\infty}$-norm of the second fundamental form such that if $$\frac{n}{V(M)^{1/p}}||H||_{2p}^2-C<\lambda_1(M)$$ then $M$ is diffeomorphic to an $n$-dimensional sphere.\\
\indent
We gave a second extrinsic pinching result in \cite{Roth}. We proved the pinching result associated with Inequality (\ref{lowerboundh}), that is, there exists a constant $C$ depending on the $L^{\infty}$-norm of the second fundamental form so that if
$$\td(R)<\frac{1}{\Hinf}+C,$$
then $M$ is diffeomorphic and almost isometric to a geodesic sphere of radius $R$.\\
\indent
In this paper, we extend the results in \cite{Roth} to high order mean curvatures $H_k$, which are the natural generalization of the mean curvature $H$. They are defined to be the $k$-th elementary sym\-metric polynomial in the principal curvature of $M$ (see Section \ref{preliminaries}). For instance, for hypersurfaces of $\R^{n+1}$, up to a multiplicative constant, $H_1$ is the mean curvature, $H_2$ is the scalar curvature and $H_n$ is the Gauss-Kronecker curvature. These curvatures give, in general, better inequalities than those involving the mean curvature $H$. Indeed, R.C. Reilly gave in \cite{Re2} a sharper upper bound for the first eigenvalue of the Laplacian for hypersurfaces of $\R^{n+1}$. The analogue of this upper bound was proved by J.F. Grosjean for hypersurfaces of $\Ss^{n+1}$ and $\HH^{n+1}$ (see \cite{Gr3}). \\
\indent
Regarding the extrinsic radius, T. Vlachos improved Inequality (\ref{lowerboundh}) in terms of $H_k$, {\it i.e.},
\beqt\label{lowerboundhk}
\td(R)^k\geqslant\frac{1}{\Hkinf},
\eeqt
with equality only for the geodesic spheres of radius $R$  (see \cite{Vla}). This is an improvement of Inequality (\ref{lowerboundh}) since we have the following sequence of inequalities
$$H_k^{1/k}\leqslant\cdots\leqslant H_2^{1/2}\leqslant H.$$
\indent
A first question is to know if Inequality (\ref{lowerboundhk}) can be improved, replacing the $L^{\infty}$-norm by an $L^p$-norm, as for (\ref{lowerboundh}) (see \cite{Roth} for details).\\
\indent 
We prove, in Section \ref{estimates}, that such $L^p$-lower bounds are true for $\delta\geqslant0$ and any $k\in\{1,\cdots,n\}$ if $H_k$ is a positive function. Moreover, the equality is characterized by geodesic hyperspheres. Then, another natural question, is to know if there exists a pinching result for these inequalities. That is, is there a constant $C$ such that if the pinching condition 
 \beqt
\td(R)^k<\frac{1}{\Hkinf}+C\tag{$P_C$}
\eeqt
holds, then $M$ is closed to a geodesic hypersphere?\\
\indent
In what follows, we denote by $\Mndrk$ {\it the family of all compact, connected and oriented $n$-dimensional Riemannian manifolds without boundary  isometrically immersed into $\R^{n+1}$ if $\delta=0$ or into the open hemisphere of $\Ss^{n+1}(\delta)$ if $\delta>0$, of extrinsic radius $R$, volume equal to $1$ and positive $H_k$}. We prove the following
\begin{thm}\label{thm1}
Let $(M^n,g)\in\Mndrk$, $\delta\geqslant0$ and $p_0$ be the center of the ball of radius $R$ containing $M$. Then for any $\varepsilon>0$, there exists a positive constant $C_{\varepsilon}$ depending only on $n$, $\delta$, the $L^{\infty}$-norm of the mean curvature and the $L^{2p}$-norm of the $k$-th mean curvature $H_k$ such that if 
\beqt
\td(R)^k<\frac{1}{\Hkp}+C_{\varepsilon}\tag{$P_{C_{\varepsilon}}$}
\eeqt
then 
\begin{enumerate}[i)]
\item $\phi(M)\subset B\big(p_0,R\big)\setminus B\big(p_0,R-\varepsilon\big).$
\item $\forall x\in S\big(p_0,R\big),\quad B(x,\varepsilon)\cap\phi(M)\neq\emptyset.$
\end{enumerate}
\end{thm}
\noindent
{\bfseries Remarks}\be[1)]
\item We will see in the proof that $ C_{\varepsilon}\lgra0$ when $\Hinf\lgra\infty$ or $\varepsilon\lgra0$.
\item An immediate consequence of $i)$ and $ii)$ of Theorem \ref{thm1} is that the Haussdorff-distance between $M$ and $S(p_0,R)$ satisfies
$$d_H\left(M,S(p_0,R)\right)\leqslant\varepsilon,$$
\ee 

If the pinching condition is strong enough, with a control on the $L^{\infty}$-norm of the second fundamental form $B$ instead of the $L^{\infty}$-norm of the mean curvature, we obtain that $M$ is diffeomorphic and almost isometric to a geodesic sphere in the following sense:
\begin{thm}\label{thm2}
Let $(M^n,g)\in\Mndrk$, $\delta\geqslant0$ and $p_0$ be the center of the ball of radius $R$ containing $M$. Then there exists a constant $C$ depending only on $n$, $\delta$, the $L^{\infty}$-norm of the second fundamental form and the $L^{2p}$-norm of the $k$-th mean curvature $H_k$ such that if
\beqt
\td(R)^k<\frac{1}{\Hkp}+C\tag{$P_C$}
\eeqt
then $M$ is diffeomorphic to $S(p_0,R)$.\\
\indent
More precisely, there exists a diffeomorphism $F$ from $M$ into the geodesic hypersphere $S(p_0,R)$ which is a quasi-isometry. That is, for all $\theta\in]0,1[$, there exist a constant $C_{\theta}$ depending on $n$, $\delta$, $\Binf$, $\Hkp$ and  $\theta$ such that the pinching condition $(P_C)$ implies
$$\big||dF_x(u)|^2-1\big|\leqslant\theta,$$
for all unit vector $u\in T_xM$.
\end{thm}
\noindent
{\bfseries Remarks.} \be[1)]
\item
It is obvious that the pinching condition 
\beqt
\td(R)^k<\frac{1}{\Hkinf}+C\tag{$\widetilde{P_C}$}
\eeqt
implies $(P_C)$. So we deduce immediately from Theorems \ref{thm1} and \ref{thm2} the same results with the second pinching condition $(\widetilde{P_C})$.
\item
In general, the constants $C_{\varepsilon}$ and  $C_{\theta}$ of Theorems \ref{thm1} and \ref{thm2} depend on $\Hkp$. In fact, the constant $C_{\varepsilon}$ of Theorem \ref{thm1} does not depend on $\Hkp$ when
\be
\item $\delta>0$
\item $\delta=0$ and $k\geqslant4$,
\item $\delta=0$ and $p\geqslant\frac{n}{2k}$,
\ee
and the constant $C_{\theta}$ of Theorem \ref{thm2} does not depend on $\Hkp$ when $\delta=0$ and $p\geqslant\frac{n}{2k}$.
\item Our approach does not work in the case $\delta<0$. As we will see, it is due to the fact that the function $\cd$ is increasing if $\delta<0$.
\item By homothety, we can deduce the same results for manifolds with arbitrary volume. Indeed, $(M,g')\in\mathcal{M}(n,\delta',R')$, with $g'=V(M)^{-2/n}g$, $\delta'= V(M)^{2/n}\delta$ and $R'=V(M)^{-1/n}R$.
\ee
These results are of special interest for $k=2$. Indeed, in that case, up to a constant, $H_2$ is the scalar curvature. So we obtain a relation between the extrinsic radius and the scalar curvature, which is an intrinsic geometric invariant. In particular, we have the following corollary
\begin{cor}\label{cor1}
Let $(M^n,g)\in\mathcal{M}(n,0,R)$ and $p_0$ be the center of the ball of radius $R$ containing $M$. We assume that $(M^n,g)$ has positive scalar curvature. Then for any $p\geqslant\frac{n}{4}$, there exists a constant $C$ depending only on $n$ and $\Binf$ such that if
$$
R^2<\frac{n(n-1)}{||\mathrm{Scal}||_{2p}}+C
$$
then $M$ is diffeomorphic and almost isometric to $S(p_0,R)$ in the sense of Theorem \ref{thm2}\\
\end{cor}
\noindent
{\bfseries Acknowledgement.} The author would like to express his gratefulness to his advisors Jean-Fran\c{c}ois Grosjean and Oussama Hijazi for their encouragement.
%%%%%%%%%%%%%%%%%%%%%%%%%%%%%%%%%%%%%%%%%%%%%%%%%%%%%%%%%%%%%%%%%%%%%%%%
%%%%%%%%%%%%%%%%%%%%%%%%%%%%%%%%%%%%%%%%%%%%%%%%%%%%%%%%%%%%%%%%%%%%%%%%
\section{Preliminraries}
\label{preliminaries}
First, let's introduce the following functions:
$$\sd(t)=\left\lbrace \begin{array}{ll}
\frac{1}{\sqrt{\delta}}\sin( \sqrt{\delta}\,t) & \text{if}\ \delta>0\\
t & \text{if}\ \delta=0 
\end{array}\right. $$
and 
$$\cd(t)=\left\lbrace \begin{array}{ll}
\cos(\sqrt{\delta}\,t) & \text{if}\ \delta>0\\
1 & \text{if}\ \delta=0 
\end{array}\right. $$
\indent
Throughout this paper, we consider a manifold $(M^n,g)\in\Mndrk$, $\delta\geqslant0$. For simplicity, we assume that $\delta=0$ or $1$. By homothety, we can deduce the results for any $\delta\geqslant0$. Let $\nu$ be the outward unit vector field. The second fundamental form $B$ of the immersion is defined by
$$B(X,Y)=\left\langle \nablab_X\nu,Y\right\rangle ,$$
where $<\cdot,\cdot>$ and $\nablab$ are respectively the Riemannian metric and the Riemannian connection of $\Md$. The mean curvature of the immersion is $$H=\frac{1}{n}\trace(B).$$\\
\indent
Now let's recall the definition of the high order mean curvature $H_k$. Let $\{e_1,\cdots,e_n\}$ be an orthonormal frame of $T_xM$. For all 
$k\in\{1,\cdots,n\}$, the $k$-th mean curvature of the immersion is
$$H_k=\left(\begin{array}{c}n\\k\end{array}\right)^{-1}\sum_{\begin{array}{c}1\leqslant i_1,\cdots,
    i_k\leqslant n\\1\leqslant j_1,\cdots,
    j_k\leqslant n\end{array}}\epsilon\left(\begin{array}{c}i_1\cdots
    i_k\\j_1\cdots
    j_k\end{array}\right)B_{i_1j_1}\cdots B_{i_kj_k},$$
    where the $B_{ij}$ are the coefficients of the real second fundamental form. The symbols $\epsilon\left(\begin{array}{c}i_1\cdots
    i_k\\j_1\cdots
    j_k\end{array}\right)$ are the usual premutation symbols which are
zero if the sets $\{i_1,\cdots,i_k\}$ and $\{j_1,\cdots,j_k\}$ are
    different or if there exist distinct $p$ and $q$ with $i_p=i_q$. For all other cases, $\epsilon\left(\begin{array}{c}i_1\cdots
    i_k\\j_1\cdots
    j_k\end{array}\right)$ is the signature of the permutation $\left(\begin{array}{c}i_1\cdots
    i_k\\j_1\cdots
    j_k\end{array}\right)$. By convention, we set $H_0=1$ et $H_{n+1}=0$.

For $k\in\{1,\cdots,n\}$, the symmetric $(1,1)$-tensor associated to $H_k$ is
$$T_k=\frac{1}{k!}\sum_{\begin{array}{c}1\leqslant i,i_1,\cdots,
    i_k\leqslant n\\1\leqslant j,j_1,\cdots,
    j_k\leqslant n\end{array}}\epsilon\left(\begin{array}{c}i_1\cdots
    i_k\\j_1\cdots
    j_k\end{array}\right)B_{i_1j_1}\cdots B_{i_kj_k}e_i^*\otimes e_j^*.$$
    This tensor is divergence free, symmetric $(1,1)$-tensor. For any symmetric $(1,1)$-tensor, we define the following function
\beqt\label{HT}
H_T(x)=\sum_{i=1}^nB_x(Te_i,e_i),
\eeqt
where $\base$ is an orthonormal frame of $T_xM$. Then, we have the following relations 
    \begin{lem}\label{trace} 
For $k\in\{1,\cdots,n\}$, we have:
\begin{enumerate}
\item $\trace(T_k)=m(k)H_k,$
\item $H_{T_k}=m(k)H_{k+1},$
\end{enumerate}
where $m(k)=(n-k)\left(\begin{array}{c}n\\k\end{array}\right)$ and $H_{T_k}$ is given by (\ref{HT}).
\end{lem}
\indent
Let $p_0\in\Md$. We denote by $r(x)=d(p_0,x)$ the geodesic distance from $p_0$ to $x$ on $(\Md,g_{can})$. We denote by $\nabla$ ({\it resp.} $\nablab$) the gradient associated with $(M,g)$ ({\it resp.} $(\Md,g_{can}))$. Let $Z:=\sd(r)\nablab r$ be the postion vector field and $Z^T=\sd(r)\nabla r$ its tangential projection on the tangent bundle of $\phi(M)$.
We have the following lemma (see \cite{Hei} or \cite{Gr2} for a proof):
\begin{lem}\label{lem1}
Let $T$ be a positive definite $(1,1)$-tensor with $\ddiv(T)=0$. Then we have
$$\displaystyle\ddiv(TZ^T)\geqslant\cd(r)\trace(T)-\left\langle Z,H_T\right\rangle.$$
If $\delta=0$ and $T$ is the identity, then equality holds.
\end{lem}
\noindent
{\bfseries Remark.} 
Note that, by integration, in the case $\delta=0$ and $T=\iid$, this lemma is nothing else but the Hsiung-Minkowski formula (see \cite{Hsi}).
\\\\
\indent
Finally, we recall the following lemma (see \cite{CG} or \cite{Roth})
\begin{lem}\label{fondlem}
Let $(M^n,g)$ be a compact, connected, oriented $n$-dimensional Riemannian manifold without boundary isometrically immersed by $\phi$ into $\R^{n+1}$ or an open hemisphere of $\Ss^{n+1}(\delta)$. Let $\xi$ be a nonnegative continuous function on $M$ such that $\xi^k$ is smooth for $k\geqslant2$. Let $0\leqslant l<m\leqslant2$ such that
$$\frac{1}{2}\xi^{2k-2}\Delta\xi^2\leqslant\ddiv\omega+(\alpha_1+k\alpha_2)\xi^{2k-l}+(\beta_1+k\beta_2)\xi^{2k-m},$$
where $\omega$ is a 1-form and $\alpha_1$, $\alpha_2$, $\beta_1$ and $\beta_2$ some nonnegative constants. Then for all $\eta>0$, there exists a constant $L$ depending only on $\alpha_1$, $\alpha_2$, $\beta_1$, $\beta_2$, $\Hinf$ and $\eta$ such that if $||\xi||_{\infty}>\eta$, then
$$||\xi||_{\infty}\leqslant L||\xi||_2.$$
Moreover, $L$ is bounded when $\eta\longrightarrow\infty$ and if $\beta_1>0$, $L\longrightarrow\infty$ when $\Hinf\longrightarrow\infty$ or $\eta\longrightarrow0$.
\end{lem}

%%%%%%%%%%%%%%%%%%%%%%%%%%%%%%%%%%%%%%%%%%%%%%%%%%%%%%%%%%%%%%%%%%%%%%%%
%%%%%%%%%%%%%%%%%%%%%%%%%%%%%%%%%%%%%%%%%%%%%%%%%%%%%%%%%%%%%%%%%%%%%%%%
\section{New estimates for the extrinsic radius}
\label{estimates}
In this section, we give new lower bounds for the extrinsic radius of hypersurfaces of the Euclidean space and the open hemisphere $\Ss^{n+1}_+(\delta)$. Let $(M^n,g)$ be a compact, connected and oriented $n$-dimensional Riemannian manifold without boundary isometrically immersed by $\phi$ in $\R^{n+1}$ or $\Ss^{n+1}_+(\delta)$. An immediate consequence of Lemma \ref{lem1} is the following
\begin{prop}\label{prop0}
Let $k\in\{1,\cdots,n\}$, if $H_k$ is positive, then for all $j\in\{1,\cdots,k\}$, the function $H_j$ is positive and 
$$\int_MH_{j-1}\,\cd(r)dv_g\leqslant\int_MH_j\,\sd(r)dv_g.$$
\end{prop}
\noindent {\it Proof:} By Lemma \ref{lem1}, we have
\beQ
\int_M\cd(r)\,\trace(T)dv_g&\leqslant&\int_M\big|\left\langle Z,H_T\right\rangle\big|dv_g\\
&\leqslant&\int_M|H_T|\,\sd(r)dv_g
\eeQ
If $T=T_{j-1}$ is the $(1,1)$-tensor associated with $H_{j-1}$, we have
\beQ
\int_M\cd(r)\,H_{j-1}dv_g&\leqslant&\int_M|H_j|\,\sd(r)dv_g
\eeQ
Moreover, Barbosa and Colares proved in \cite{BC} that if $H_k$ is positive, then for all $j\in\{1,\cdots,k\}$, $H_j$ is positive. So for all $j\in\{1,\cdots,k\}$, we have
$$\int_MH_{j-1}\,\cd(r)dv_g\leqslant\int_MH_j\,\sd(r)dv_g.$$
$\finpreuve$\\
From this proposition, we deduce the following estimates
\begin{thm}\label{thm5}
Let $(M^n,g)\in\Mndrk$ (we do not assume that $V(M)=1$). Then for all $p\geqslant1$, we have
\beqt\label{lowerboundhkp}
\td(R)^k\geqslant\dfrac{V(M)^{1/p}}{||H_k||_p}.
\eeqt
Equality occurs if and only if $M$ is a geodesic hypersphere.
\end{thm}
\noindent
{\bfseries Remarks.}
\be[i)]
\item These lower bounds improve Inequalities $(\ref{lowerboundh})$ and $(\ref{lowerboundhk})$.
\item In the case $k=2$, this inequality \ref{lowerboundhkp} translates to
\beqt\label{Deshmukh}
R^2 \geqslant \dfrac{n(n-1)}{||\sscal||_1},\ \text{if}\ \delta=0
\eeqt
\beqt\label{Aliasbetter}
\tan^2(R)\geqslant \dfrac{n(n-1)}{||\sscal-n(n-1)||_1},\ \text{if}\ \delta=1.
\eeqt
Inequaliy (\ref{Deshmukh}) was proved by S. Deshmukh (\cite{Des}). In the spherical case, our inequality improves the following one due to L.J Alias (\cite{Ali})
\beqt\label{Alias}
\sin^2(R)\geqslant\dfrac{n(n-1)}{||\sscal||_1}.
\eeqt
Indeed, Alias proved Inequality (\ref{Alias}) with the assumption that \\$\Ric(M)\geqslant (n+2)(n-1)$. With such an assumption, Inequalities (\ref{Alias}) and (\ref{Aliasbetter}) are exactly the same, but Inequality (\ref{Aliasbetter}) is valid without any assumption on the Ricci curvature of $M$.
\ee
\noindent {\it Proof:} We use Proposition \ref{prop0} and the fact that the functions $\sd$ and $\cd$ are respectively increasing and decreasing.
\beQ
\cd(R)^kV(M)&\leqslant&\cd(R)^{k-1}\int_M\cd(r)dv_g\\
&\leqslant&\cd(R)^{k-1}\int_MH\,\sd(r)dv_g\\
&\leqslant&\cd(R)^{k-2}\sd(R)\int_MH\,\cd(r)dv_g\\
&\leqslant&\cd(R)^{k-2}\sd(R)\int_MH_2\,\sd(r)dv_g\\
&\leqslant&\cdots\\
&\leqslant&\sd(R)^k\int_MH_kdv_g.
\eeQ
So we have
$$\td(R)^k\geqslant\frac{V(M)}{||H_k||_1},$$
and equality occurs if and only if $r=R$ for all $x\in M$.
Finally, the H\"older inequality gives the result with the $L^p$-norm, for all $p\geqslant1$.
$\finpreuve$
%%%%%%%%%%%%%%%%%%%%%%%%%%%%%%%%%%%%%%%%%%%%%%%%%%%%%%%%%%%%%%%%%%%%%%%%

\section{An $L^2$-approach to pinching}
Let $(M^n,g)\in\Mndrk$. A first step in the proof of the pinching restults is to prove that the pinching condition
\beqt
\td(R)^k<\frac{1}{\Hkp}+C\tag{$P_C$}
\eeqt
implies that $M$ is closed to a geodesic hypersphere in an $L^2$-sense.
For this, let's introduce the following functions:
$$\left\lbrace \begin{array}{l}
\varphi(x)=\sd^2(R)-\sd^2(r), \\ 
\psi(x)=\cd(r)|Z^T|.
\end{array}\right. $$
In what follows, we assume that the pinching constant satisfies $C<1$. We prove the following lemma
\begin{lem}\label{lemphi}
The pinching condition $(P_C)$ with $C<1$ implies 
$$||\varphi||_2^2\leqslant A_1C,$$
where $A_1$ is a positive explicit constant depending only on the dimension $n$, $\delta$ and $\Hkp$. Moreover in certain cases, the dependence on $\Hkp$ can be replaced by a dependence on $\Hinf$, precisely when
\be
\item $\delta>0$,
\item $\delta=0$ and $k\geqslant4$,
\item $\delta=0$ and $p\geqslant\frac{n}{2k}$. 
\ee
\end{lem}
\noindent {\it Proof:} Since $\sd$ and $\cd$ are respectively increasing and decreasing, we have
\beQ
||\varphi||_2^2&\leqslant&\sd^2(R)\int_M\left( \sd^2(R)-\sd^2(r)\right) \\
&\leqslant&\sd^2(R)\left[ \td^2(R)\left( \int_M\cd(r)\right)^2-\int_M\sd^2(r)\right]
\eeQ
By the H\"older inequality, we get
\beQ
||\varphi||_2^2&\leqslant&\sd^2(R)\left[ \td^2(R)\left( \int_M\cd(r)\right)^2-\frac{1}{\HHkp}\left( \int_MH_k\sd(r)\right)^2\right]  
\eeQ
Now, using Proposition \ref{prop0}, we have
\beQ
||\varphi||_2^2&\leqslant&\sd^2(R)\left[ \td^2(R)\left( \int_M\cd(r)\right)^2-\frac{1}{\HHkp}\left( \int_MH_{k-1}\cd(r)\right)^2\right]\\
&\leqslant&\sd^2(R)\left[ \td^2(R)\left( \int_M\cd(r)\right)^2-\frac{1}{\HHkp\td^2(R)}\left( \int_MH_{k-1}\sd(r)\right)^2\right]\\\\
\eeQ
Proposition \ref{prop0} applied $(k-2)$ more times, yields
\beQ
||\varphi||_2^2&\leqslant&\sd^2(R)\left[ \td^2(R)\left( \int_M\cd(r)\right)^2-\frac{1}{\HHkp\td^{2k-2}(R)}\left( \int_M\cd(r)\right)^2\right]\\
&\leqslant&\frac{\sd^2(R)}{\td^{2k-2}(R)}\left( \int_M\cd(r)\right)^2\left[\td^{2k}(R)- \frac{1}{\HHkp}\right] 
\eeQ
Since we assume $(P_C)$ is true with $C<1$, if $\delta>0$, we have
\beQ
||\varphi||_2^2&\leqslant&\frac{C}{\delta\td^{2k-2}(R)}\left[ 1+\frac{2}{\Hkp}\right] 
\eeQ
Moreover, we have $$\frac{1}{\td(R)}\leqslant\Hkp^{1/k},$$
so
\beQ
||\varphi||_2^2&\leqslant&\frac{C}{\delta}\left[ \Hkp^{\frac{2k-2}{k}}+2\Hkp^{\frac{k-2}{k}}\right]
\eeQ
Since $H_k>0$, we have $H\geqslant H_k^{1/k}$ and then $\Hkp\leqslant\Hinf^{k}$. Finally, for $k\geqslant2$, we have
\beQ
||\varphi||_2^2&\leqslant&\frac{C}{\delta}\left[ \Hinf^{2k-2}+2\Hinf^{k-2}\right]:=A_1C,
\eeQ
If $\delta=0$, then
\beQ
||\varphi||_2^2&\leqslant&\frac{C}{R^{2k-4}}\left[ 1+\frac{2}{\Hkp}\right]\\
&\leqslant&C\left[ \Hkp^{\frac{2k-4}{k}}+2\Hkp^{\frac{k-4}{k}}\right]:=A_1C,
\eeQ
If, in addition, $4\leqslant k\leqslant n$, by the same argument as above,  $A_1$ depends on $n$ and $\Hinf$.\\
Moreover, for $\delta=0$, we have the following lower bound for the $k$-th mean curvature (see \cite{BZ} page 221)
$$\int_MH_k^{n/k}dv_g\geqslant\omega_n,$$
where $\omega_n$ is the volume of the $n$-dimensional Euclidean sphere. Then, for $p\geqslant\frac{n}{2k}$, we have
$$\omega_n^{k/n}\leqslant||H_k||_{2p}\leqslant\Hinf^{k}.$$
So the dependence on $\Hkp$ can be replaced by a dependence on $\Hinf$.
$\finpreuve$
\begin{lem}\label{lempsi}
The pinching condition $(P_C)$ implies 
$$||\psi||_2^2\leqslant A_2C,$$
where $A_2$ is a positive explicit constant depending only on $n$ and $\Hinf$.
\end{lem}
\noindent {\it Proof:}
Since $\cd(r)\leqslant1$, we have
\beQ
||\psi||_2^2&\leqslant&\int_M|Z^T|^2\leqslant\int_M|Z|^2-\left\langle Z,\nu\right\rangle^2\\
&\leqslant&\int_M\sd^2(r)-\frac{1}{\HHkp}\left( \int_MH_k\left\langle Z,\nu\right\rangle \right)^2
\eeQ
By Lemma \ref{lem1}, we have
\beQ
||\psi||_2^2&\leqslant&\td^2(R)\left( \int_M\cd(r)\right)^2-\frac{1}{\HHkp}\left( \int_M\cd(r)H_{k-1}\right) ^2\\
&\leqslant&\td^2(R)\left( \int_M\cd(r)\right)^2-\frac{1}{\td^2(R)\HHkp}\left( \int_M\sd(r)H_{k-1}\right)^2
\eeQ
Using Proposition \ref{prop0}, we get
\beQ
||\psi||_2^2&\leqslant&\td^2(R)\left( \int_M\cd(r)\right)^2-\frac{1}{\td^2(R)\HHkp}\left( \int_M\cd(r)H_{k-2}\right)^2\\
&\leqslant&\td^2(R)\left( \int_M\cd(r)\right)^2-\frac{1}{\td^4(R)\HHkp}\left( \int_M\cd(r)H_{k-2}\right)^2\\\\
&\leqslant&\cdots\\\\
&\leqslant&\td^2(R)\left( \int_M\cd(r)\right)^2-\frac{1}{\td^{2k-2}(R)\HHkp}\left( \int_m\cd(r)\right)^2\\
&\leqslant&\frac{1}{\td^{2k-2}(R)}\left( \td^2(R)-\frac{1}{\HHkp}\right) 
\eeQ
So the pinching condition $(P_C)$ with $C<1$ implies
\beQ
||\psi||_2^2&\leqslant&\frac{C}{\td^{2k-2}(R)}\left( 1+\frac{2}{\Hkp}\right) \\
&\leqslant&C\Hkp^{\frac{2k-2}{k}}\left( 1+\frac{2}{\Hkp}\right)\\
&\leqslant&C\left( \Hkp^{\frac{2k-2}{k}}+2\Hkp^{\frac{k-2}{k}}\right) \\
&\leqslant&C\left( \Hinf^{2k-2}+2\Hinf^{k-2}\right):=A_2C,
\eeQ
$\finpreuve$\\
\indent
The next step to prove Theorems \ref{thm1} and \ref{thm2} is to get $L^{\infty}$-estimates from these $L^2$-estimates. For this, we will use Lemma \ref{fondlem}.
%%%%%%%%%%%%%%%%%%%%%%%%%%%%%%%%%%%%%%%%%%%%%%%%%%%%%%%%%%%%%%%%%%%%%%%%
\section{Proof of Therorem \ref{thm1}}
The proof of Theorem \ref{thm1} is an immediate consequence of the following three lemmas
\begin{lem}\label{prop1}
For any $\varepsilon>0$, there exists $K_{\varepsilon}$ depending on $n$ and $\Hinf$ so that if $(P_{K_{\varepsilon}})$  is true, then
$$\phi(M)\subset B(p_0,R)\setminus B(p_0,R-\varepsilon).$$
Moreover, $K_{\varepsilon}\lgra0$ when $\Hinf\lgra+\infty$ or $\varepsilon\lgra 0$.
\end{lem}
\noindent {\it Proof:} We showed in \cite{Roth} that the function $\varphi$ satisfies
\beqt
\varphi^{2k-2}\Delta\varphi^2\leqslant\ddiv(\omega)+(\alpha_1+k\alpha_2)\varphi^{2k-1}+(\beta_1+k\beta_2)\varphi^{2k-2},
\eeqt
 where $\omega$ is a 1-form, $\alpha_1$, $\alpha_2$, $\beta_1$ and $\beta_2$ some nonnegative constants depending on $n$, $\delta$ and $\Hinf$. We can apply Lemma \ref{fondlem} to the function $\varphi$ with $l=1$ and $m=2$. We deduce that if $||\varphi||_{\infty}>\varepsilon$ then there exists a constant $L$ such that $$||\varphi||_{\infty}\leqslant L||\varphi||_2.$$
On the other hand, by Lemma \ref{lemphi}, we know that if the pinching condition $(P_C)$ is satisfied for $C\leqslant1$, then 
$$||\varphi||_2^2\leqslant A_1C.$$
Take $C=K_{\varepsilon}=\iinf\left\lbrace 1,\frac{\varepsilon^2}{L^2A_1}\right\rbrace$. This choice implies
$$ ||\varphi||_{\infty}\leqslant\varepsilon,$$
that is, $\td^2(R)-\td^2(r)\leqslant\varepsilon$. Finally, we can choose $K_{\varepsilon}$ smaller in order to have
$R-r\leqslant\varepsilon.$
$\finpreuve$\\
The second lemma is due to B. Colbois and J.F. Grosjean (see \cite{CG}).
\begin{lem}\label{lem3}
Let $x_0$ be a point of the sphere $S(0,R)$ of $\R^{n+1}$. Assume that $x_0=Ru$ where $u\in\Ss^n$. Now let $(M^n,g)$ be a compact, connected and oriented $n$-dimensional Riemannian manifold without boundary isometrically immersed by $\phi$ into $\R^{n+1}$ so that $$\phi(M)\subset\Big( B(p_0,R+\eta)\setminus B(p_0,R-\eta)\Big) \setminus B(x_0,\rho)$$ with $\rho=4(2n-1)\eta$ and suppose there exists a point $p\in M$ so that $\left\langle Z,u\right\rangle(p)\geqslant0$.  Then there exists $y_0\in M$ so that the mean curvature satisfies $|H(y_0)|>\frac{1}{4n\eta}$. 
\end{lem}
\noindent
{\bfseries Remark.}
Note that in \cite{CG}, it is supposed that $\left\langle Z,u\right\rangle(p) >0$, but the condition $\left\langle Z,u\right\rangle(p)\geqslant0$ is sufficient.
\\\\
In \cite{Roth}, we give a corresponding lemma for the spherical case.
\begin{lem}\label{lem4}
Let $x_0$ be a point of the sphere $S(p_0,R)$ of an open hemisphere of $\Ss^{n+1}(\delta)$. Let $(M^n,g)$ be a compact, connected and oriented $n$-dimensional Riemannian manifold without boundary isometrically immersed by $\phi$ into this open hemisphere of $\Ss^{n+1}(\delta)$ so that $$\phi(M)\subset\Big( B(p_0,R)\setminus B(p_0,R-\eta)\Big) \setminus B(x_0,\rho)$$
with $\rho$ such that 
$$ \td(R/2)-\td\left( (R-\rho)/2\right) =4(2n-1)\eta$$
 Then there exists two constants $D$ and $E$ depending on $n$, $\delta$ and $R$ such that if $\eta\leqslant D$, then there exists $y_0\in M$ so that
$$| H(y_0)|\geqslant \frac{E}{8n\eta}.$$
\end{lem}
\noindent
{\bfseries Proof of Theorem \ref{thm1}:}\\
The proof for $\delta=0$ is an immediate consequence of Lemmas \ref{prop1} and \ref{lem3}. Indeed, let $\varepsilon>0$. By Lemma \ref{prop1}, there exists $K_{\varepsilon}$ such that if $(P_{K_{\varepsilon}})$ is true, then
$$\phi(M)\subset B(p_0,R)\setminus B(p_0,R-\varepsilon).$$
Let $x=Ru\in S(0,R)$ and assume that $\phi(M)\cap B(x,\varepsilon)=\emptyset$. Since $R$ is the extrinsic radius, there exists a point $p\in M$ such that $\left\langle Z,u\right\rangle(p)\geqslant0$. If $(P_{C_{\varepsilon}})$ is true with $C_{\varepsilon}=K_{\frac{\varepsilon}{4(2n-1)}}$ and $\varepsilon<\frac{2}{3\Hinf}$, then by Lemma \ref{lem3}, there exists $y_0\in M$ so that
$$|H(y_0)|\geqslant\frac{1}{4n(\frac{\varepsilon}{4(2n-1)})}\geqslant\frac{2n-1}{n\varepsilon}>\Hinf,$$
which is a contradiction. Finally, $\phi(M)\cap B(x,\varepsilon)\neq\emptyset$ which completes the proof for $\delta=0$.\\
For $\delta\neq0$, let $\varepsilon>0$. We set $0<\eta:=\iinf\left\lbrace D,\varepsilon,\frac{\gamma(\varepsilon)}{8(2n-1)}\right\rbrace$, where
$$\gamma(\varepsilon)=\td\left( \frac{R}{2}\right)-\td\left( \frac{R-\varepsilon}{2}\right).$$
Note that $\gamma$ is an increasing smooth function with $\gamma(0)=0$. From Lemma \ref{prop1}, there exists $C_{\varepsilon}=K_{\eta}$ such that $(P_{C_{\varepsilon}})$ implies $$R-r\leqslant\eta\leqslant\varepsilon.$$ That's the first point of Theorem \ref{thm1}. Assume that $\varepsilon<\gamma^{-1}\left( \frac{2E}{3\Hinf}\right)$. Suppose there exists $x\in S(p_0,R)$ such that $B(x,\varepsilon)\cap M=\emptyset$. Since $\gamma(\varepsilon)\geqslant 4(2n-1)\eta$, by Lemma \ref{lem4}, there exists a point $y_0\in M$ so that
$$|H(y_0)|\geqslant\frac{E}{8n\eta}\geqslant\frac{(2n-1)E}{n\gamma(\varepsilon)}>\Hinf.$$
Hence a contradiction and $B(x,\varepsilon)\cap M\neq\emptyset$. \\
Moreover, for $\delta=0$ or $\delta>0$, by Lemma \ref{prop1}, $C_{\varepsilon}\lgra0$ when \\
$\Hinf\lgra0$ or $\varepsilon\lgra0$.
$\finpreuve$
\section{Proof of Theorem \ref{thm2}}
We first need the following lemma
\begin{lem}\label{lem2}
 For any $\varepsilon>0$, there exists $K_{\varepsilon}$ depending on $n$ and $\Binf$ so that if $(P_{K_{\varepsilon}})$  is true, then
 $$||\psi||_{\infty}\leqslant\varepsilon.$$
 Moreover, $K_{\varepsilon}\lgra0$ when $\Binf\lgra+\infty$ or $\varepsilon\lgra 0$.
\end{lem}
\noindent {\it Proof:} We proved in \cite{Roth} that the function $\psi$ satisfies
\beqt
\psi^{2k-2}\Delta\psi^2\leqslant\ddiv(\omega)+(\alpha_3+k\alpha_4)\psi^{2k-1}+(\beta_3+k\beta_4)\psi^{2k-2}.
\eeqt
Now applying Lemma \ref{fondlem} with $l=1$ and $m=2$, we get that for any $\eta>0$, there exists $L$ depending on $n$ and $\Binf$ so that if $||\psi||_{\infty}>\eta$ then 
$$||\psi||_{\infty}\leqslant L||\psi||_2.$$
From Lemma \ref{lemphi}, we know that if $(P_C)$ holds, then
$$||\psi||_2^2\leqslant A_2C.$$
Let $\varepsilon>0$, we set $C=K_{\varepsilon}=\iinf\left\lbrace 1,\frac{\varepsilon^2}{L^2A_2}\right\rbrace$. For this choice of $C$ we get 
$$||\psi||_{\infty}\leqslant\varepsilon.$$
This completes the proof.
$\finpreuve$\\
{\bfseries Proof of Theorem\ref{thm2}:}\\
We consider the following map
\function{F}{M}{S(p_0,R)}{x}{\expo\left(R\left( \diff\expo\right) ^{-1}(\nablab r)\right) .}
We proved in \cite{Roth} that for any uniatry vector $u\in T_xM$
 \beQ
 \Big|\big|\diff_xF(u)\big|^2-1\Big|&\leqslant&\frac{1}{\sd^2(r)}\big|\sd^2(R)-\sd^2(r)\big|+\frac{\sd^2(R)}{\cd(r)\sd^3(r)}||\psi||_{\infty}.
\eeQ
From Lemma \ref{lem3}, we know that for any $\eta>0$, there exists a constant $K_{\eta}$ so that $(P_{K_\eta})$ implies $||\psi||_{\infty}\leqslant\eta$. Moreover, by Theorem \ref{thm1}, there exist $C_{\eta}$ depending on $n$, $\delta$, $\Hinf$ and $\eta$ so that $(P_{C_{\eta}})$ implies $R-r\leqslant\eta$. We set $C'_{\eta}=\iinf\{C_{\eta},K_{\eta}\}$. Then, since $R$ is bounded by a constant depending only on $n$, $\delta$, $\Hkp$ and $\Hinf$ there exist three positive constants $A_3$, $A_4$ and $A_5$ depending on $n$, $\delta$, $\Hkp$ and $\Hinf$ so that 
 \beQ
 \Big|\big|\diff_xF(u)\big|^2-1\Big|&\leqslant&A_3||R-r||_{\infty}+A_4||\psi||_{\infty}\\
 &\leqslant&A_3\eta+A_4\eta\leqslant A_5\eta
 \eeQ
 Now, choosing $\eta=\frac{\theta}{A_5}$, we get
 \beqt\label{df6}\Big|\big|\diff_xF(u)\big|^2-1\Big|\leqslant\theta.\eeqt
For $\theta\in]0,1[$, by (\ref{df6}), $F$ is a local diffeomorphism from $M$ to $S(p_0,R)$. Since for $n\geqslant2$, $S(p_0,R)$ is simply connected, $F$ is a diffeomorphism. Moreover, the relation (\ref{df6}) says that $F$ is a quasi-isometry.$\finpreuve$\\
{\bfseries Remark.}
If $\delta=0$ and $p\geqslant\frac{n}{2k}$, we saw that
$$\omega_n^{k/n}\leqslant||H_k||_{2p}\leqslant\Hinf^{k}\leqslant\frac{1}{n^{k/2}}\Binf^k,$$
so there is no dependence on $\Hkp$.

\bibliographystyle{amsplain}
\bibliography{extrinsic2} 
\end{document}